\documentclass[a4paper,11pt]{article}
\usepackage{latexsym}
\usepackage{amssymb}
\usepackage{amsfonts}
\usepackage{amsmath}
\usepackage{indentfirst}
\usepackage{graphicx}
\usepackage{epsfig}
\usepackage{slashbox,colortbl}
\usepackage{color}
\usepackage{fullpage}

\newtheorem{prop}{Proposition}[section]
\newtheorem{teor}{Theorem}[section]

\newcommand{\cvd}{\hfill $\blacksquare$\bigskip}
\date{}

\newcommand\ptau{\ensuremath{{'\kern-.4ex\tau}}}
\newcommand\psigma{\ensuremath{{'\kern-.4ex\sigma}}}
\newcommand\ptp{\ensuremath{{'\kern-.4ex\tau'}}}

\newcommand\pC{\ensuremath{{'\kern-.2ex C}}}

\newcommand\tauone{\ensuremath{{\tau^{(1)}}}}
\newcommand\tautwo{\ensuremath{{\tau^{(2)}}}}
\newcommand\tauonetwo{\ensuremath{{\tau^{(1,2)}}}}
\newcommand\Cone{\ensuremath{{C^{(1)}}}}
\newcommand\Ctwo{\ensuremath{{C^{(2)}}}}

\newcommand\ist{\ensuremath{[\sigma,\tau]}}

\newcommand\mst{\ensuremath{\mu(\sigma,\tau)}}

\def\m#1,#2,{\ensuremath{\mu(#1,#2)}}
\def\mb#1,#2,{\ensuremath{\bar\mu(#1,#2)}}

\author{Antonio Bernini\thanks{Dipartimento di Matematica e Informatica ``U. Dini", viale Morgagni 65, University
of Firenze, Firenze, Italy, {\tt antonio.bernini@unifi.it,
luca.ferrari@unifi.it}. Partially supported by MIUR PRIN 2010-2011
grant ``Automi e Linguaggi Formali: Aspetti Matematici e
Applicativi", code H41J12000190001.}\and Luca Ferrari$^{\ast
,}$\thanks{Partially supported by INdAM-GNCS 2014 project ``Studio
di pattern in strutture combinatorie".}}

\title{Vincular pattern posets and the M\"obius function of the quasi-consecutive pattern poset}

\begin{document}

\maketitle

\begin{abstract}
We introduce vincular pattern posets, then we consider in
particular the quasi-consecutive pattern poset, which is defined
by declaring $\sigma \leq \tau$ whenever the permutation $\tau$
contains an occurrence of the permutation $\sigma$ in which all
the entries are adjacent in $\tau$ except at most the first and
the second. We investigate the M\"obius function of the
quasi-consecutive pattern poset and we completely determine it for
those intervals $\ist$ such that $\sigma$ occurs precisely once in
$\tau$.
\end{abstract}

\section{Introduction}

The study of patterns in permutations is one of the most active
trends of research in combinatorics. The richness of the notion of
permutation patterns is especially evident from its plentiful
appearances in several very different disciplines, such as
algebra, geometry, analysis, theoretical computer science, and
many others. Even if it is arguably not possible to encompass all
possible applications of this notion into a simple formal
environment, it seems reasonable to assert that the single
mathematical structure which best catches the concept of a pattern
and allows us to express a great deal of results about it is the
\emph{permutation pattern poset}.

Given two permutations $\sigma ,\tau$, we say that $\sigma \leq
\tau$ in the permutation pattern poset whenever there is an
occurrence of $\sigma$ in $\tau$ as a classical pattern. An
extremely challenging open problem concerning the permutation
pattern poset is the determination of its M\"obius function. The
problem, originally posed by Wilf \cite{W}, received quite
recently much attention, and some partial results have been
achieved \cite{SV,ST,BJJS,MS,Sm1,Sm2}. However, a complete description of
the M\"obius function of the permutation pattern poset is not yet
available. The same problem can be formulated for the
\emph{consecutive pattern poset}, where, by definition, $\sigma
\leq \tau$ whenever $\sigma$ appears in $\tau$ as a
\emph{consecutive} pattern. This poset is much easier than the
classical one and, in particular, its M\"obius function is now
completely understood \cite{BFS,SW}.

We recall here that the M\"obius function $\mu$ of a poset is an
important element of the incidence algebra of that poset, i.e. a
function mapping an interval of the poset into a scalar. More
specifically, it is possible to give a recursive definition of
$\mu$ as follows:
\begin{eqnarray*}
\mu (x,x)&=&1 \\
\mu (x,y)&=&-\sum_{x\leq z<y}\mu (x,z),\quad \textnormal{ if
$x\neq y$}.
\end{eqnarray*}

In particular, if $x\nleq y$, then $\mu (x,y)=0$. Using a duality
argument, it is possible to show that the M\"obius function of an
interval can be equally computed ``from top to bottom", according
to the following formula (for $x\neq y$):
$$
\mu(x,y)= -\sum_{x<z\le y}{\mu(z,y)}.
$$

Both formulas for computing $\mu$ will be frequently used
throughout the paper.

Consecutive and classical patterns are special (actually,
extremal) cases of the more general notion of \emph{vincular
patterns}. An occurrence of a vincular pattern is an occurrence of
that pattern in which entries are subject to certain adjacency
conditions (see next section for a precise definition). Vincular
patterns were introduced by Babson and Steingr\'imsson \cite{BS}
(under the name of generalized patterns), and constitute a vast
intermediate continent between the two lands of consecutive
patterns and classical patterns. Since the determination of the
M\"obius function is still open in the classical case and totally
solved in the consecutive one, it is conceivable that, if one is
able to define a reasonable poset structure on permutations
depending on the type of vincular patterns under consideration,
then the resulting class of permutation posets (which would
somehow interpolate between the consecutive and the vincular
posets) may shed light on the differences and on the analogies
between the two extremal cases.

In the present paper we propose a definition of \emph{vincular pattern
poset of type $A$} (whose elements are vincular patterns, or dashed permutations,
as defined in Section \ref{posets}), where $A$ is a suitable infinite matrix
whose $n$-th row encodes the type of the pattern of length $n$ belonging to the poset
(by specifying which entries of the pattern have to be adjacent in the corresponding
dashed permutation). This is done in Section \ref{posets}, where we also
address the problem of understanding in which
cases (that is, for which matrices $A$) the partial order thus
defined perfectly catches the notion of an $A$-vincular pattern; even
if we are not able to give a complete characterization of such
matrices, we find some partial results, and in particular we
single out a special class of matrices for which everything works
in the best possible way. In Section \ref{quasi} we consider a single case,
probably the closest to the consecutive one, in which all the
entries of a pattern are required to be adjacent, except at most
the first and the second. In Section \ref{mobius_quasi} we address the problem of the
computation of the M\"obius function in this poset, which we call
the \emph{quasi-consecutive pattern poset}, and we
completely determine the value $\mst$ when $\sigma$ occurs
precisely once in $\tau$. In spite of its closeness with the
consecutive case, it seems that, in the general case, the
computation $\mst$ in the quasi-consecutive pattern poset is
considerably more complicated. This is underlined, for instance, by
the fact (hinted at in the last section) that the absolute value
attained by the M\"obius function seems to be unbounded, whereas
in the consecutive pattern poset the only possible values are
$-1,0,1$.

\section{Vincular pattern containment orders}\label{posets}

Denote with $S$ the set of all finite permutations. Elements of
$S$ are represented in one-line notation, so that $\pi =\pi_1
\cdots \pi_n$ is the permutation of length $n$ mapping $i$ to
$\pi_i$, for all $i\leq n$. A \emph{dashed permutation} is a
permutation in which some dashes are possibly inserted between any
two consecutive elements. For instance, $5-13-42$ is a dashed
permutation (of length 5). The \emph{type} of a dashed permutation
$\pi$ of length $n$ is the $(0,1)$-vector $r=(r_1 ,\ldots
,r_{n-1})$ (having $n-1$ components) such that, for all $i\leq
n-1$, $r_i =0$ whenever there is no dash between $\pi_i$ and
$\pi_{i+1}$, and $r_i =1$ whenever there is a dash between $\pi_i$
and $\pi_{i+1}$. For example, the above dashed permutation
$5-13-42$ has type $(1,0,1,0)$. We remark that, in different
sources, this notion of type is expressed using a different
formalism, namely by recording the lengths of each interval of
adjacent elements in the dashed permutation. For instance, the
above dashed permutation is said to have type $(1,2,2)$.

Let $\pi$ be a permutation of length $n$ and $\rho$ be a dashed
permutation of length $k\leq n$. We say that $\pi$ \emph{contains
an occurrence of the vincular pattern} $\rho$ when there exists a
subsequence $\pi_{i_1},\ldots ,\pi_{i_k}$ of elements of $\pi$
which is order-isomorphic to $\rho$ and such that $\pi_{i_j}$ and
$\pi_{i_{j+1}}$ appear consecutively inside $\pi$ if there is no
dash between $\rho_j$ and $\rho_{j+1}$. In this case we also say that
$\rho$ is a vincular patter of $\pi$. Vincular patterns were
introduced in \cite{BS} (where they are called \emph{generalized
patterns}).

\medskip

Let $A$ be an infinite lower triangular $(0,1)$-matrix, and denote
with $r_i$ the $i$-th row vector of $A$. Given $\pi \in S$, we say
that $\rho \in S$ is an \emph{occurrence of a vincular pattern of
type $A$} (or an \emph{occurrence of an $A$-vincular pattern}) in
$\pi$ when, given that $\rho$ has length $k$, $\pi$ contains an
occurrence of $\rho$ as a vincular pattern of type
(given by the first $k-1$ entries of) $r_{k-1}$. In
this case, we also write $\rho \in_A \pi$. If $\pi$ does not
contain any occurrence of an $A$-vincular pattern $\rho$, that is
$\rho \notin_A \pi$, we say that $\pi$ \emph{avoids} $\rho$
\emph{as an $A$-vincular pattern}. The most studied special
cases are obtained when all the entries of $A$ below and on the
main diagonal are equal to 1 and when $A$ is the null matrix. In
the former case we recover the notion of classical pattern,
whereas in the latter one we get the notion of consecutive
pattern.

\medskip

In some cases, the notion of an $A$-vincular pattern can be
described by means of a suitable partial order. Given $\pi ,\rho
\in S$, having lengths $n$ and $n-1$ respectively, we say that
$\pi$ \emph{covers} $\rho$ whenever $\rho$ appears as an
$A$-vincular pattern in $\pi$. In this case, we write $\rho
\prec_A \pi$, or simply $\rho \prec \pi$ when $A$ is clear from
the context. The transitive and reflexive closure of this covering
relation is a partial order which will be called the
\emph{A-vincular pattern containment order}, and the resulting
poset will be called the \emph{A-vincular pattern poset}. When
$\sigma$ is less than or equal to $\tau$ in the $A$-vincular
pattern poset, we will write $\sigma \leq_A \pi$ or simply $\sigma
\leq \pi$, when no confusion is likely to arise.

\medskip

It is not difficult to realize that, in the two special cases
mentioned above, the resulting posets are well known. When $A$ is
the lower triangular matrix having all 1's below and on the main
diagonal we obtain the classical pattern poset, whereas when $A$
is the null matrix we get the consecutive pattern poset. We
observe that, in these two cases, the partial order relation
$\leq_A$ coincides with the binary relation $\in_A$, that is
$\sigma \leq_A \tau$ if and only if $\sigma \in_A \tau$. However
this is not always true, for the relation $\in_A$ is not
transitive in general.
In this direction, it would be interesting to characterize all
matrices $A$ for which the partial order relation (rather than the
covering one) is directly defined in terms of occurrences of
$A$-vincular patterns. More precisely, it would be nice to have a
characterization of those matrices $A$ such that $\sigma \leq_A
\tau$ if and only if $\sigma \in_A \tau$. Such matrices $A$ are
those for which the structure of the $A$-pattern poset perfectly
describes the notion of an $A$-vincular pattern. Unfortunately, we
have not been able to find such a characterization yet, so we
leave it as our first open problem. We have however some partial
results, which we are going to illustrate below.

\medskip

The first thing we observe is that in general, in the above
conjectured equivalence, $\sigma \in_A \tau$ does not imply
$\sigma \leq_A \tau$, nor does $\sigma \leq_A \tau$ imply $\sigma
\in_A \tau$, as shown by the following two examples.

\medskip

\emph{Examples.}

\begin{enumerate}

\item Let $A$ be the (lower triangular) matrix all of whose
elements are 0, except for row 3, which is $r_3 =(0,1,0)$. If
$\sigma =1234$ and $\tau =342156$, then $\sigma \in_A \tau$ (there
is precisely one occurrence, in the subsequence 3456), however
$\sigma \nleq_A \tau$ (any $A$-vincular pattern of $\tau$ of
length 5 has to be consecutive in $\tau$). Notice that this
example can be easily generalized to a matrix $A$ whose $n$-th row
has only one 1 and all successive rows are identically 0.

\item If $A$ is such that $r_4 =(1,0,0,0)$ and $r_5 =(0,0,0,0,1)$,
then it is immediate to verify that $31524\prec_A 361524$ and
$361524\prec_A 3615274$; thus we have $31524\leq_A 3615274$, but
$31524\notin_A 3615274$, since the unique occurrence is of type
$(1,0,0,1)$.

\end{enumerate}

If we restrict to the class of infinite lower triangular matrices
having constant columns, we can show that one of the previous
implications holds. Observe that, in this case, we can completely
describe $A$ using the (infinite) vector $a$ whose $i$-th
component $a_i$ is the unique nonzero value appearing in column $i$ of $A$. In
particular, we modify the notations accordingly, by writing
$\sigma \in_a \tau$ and $\sigma \leq_a \tau$ in place of $\sigma
\in_A \tau$ and $\sigma \leq_A \tau$, respectively.

\begin{prop} If $\sigma \in_a \tau$, then $\sigma \leq_a \tau$.
\end{prop}

\emph{Proof.}\quad Fix an occurrence of $\sigma$ in $\tau$ as an
$a$-vincular pattern. Let $\pi$ be the smallest consecutive
pattern of $\tau$ containing such an occurrence. Then $\pi \leq_a
\tau$: this follows from the more general fact that, given any
matrix $A$, for any consecutive pattern $\rho'$ of a permutation
$\rho$, we have $\rho' \leq_A \rho$ (starting from $\rho$,
repeatedly remove either the first or the last element until
getting to $\rho'$). To conclude the proof we now need to find a
(descending) chain of coverings from $\pi$ to $\sigma$ in the
$a$-vincular pattern poset. To this aim, partition the elements of
the fixed occurrence of $\sigma$ in $\pi$ into blocks of
consecutive elements of $\pi$ of maximal length, and suppose that
such lengths are $\alpha_1 ,\alpha_2 ,\ldots ,\alpha_r$. This
implies that the vector $a$ certainly has a 1 in positions
$\alpha_1 ,\alpha_1 +\alpha_2 ,\ldots ,\alpha_1 +\alpha_2 +\cdots
+\alpha_r$. Therefore we can remove the leftmost element of $\pi$
not belonging to the selected occurrence of $\sigma$, which is the
$(\alpha_1 +1)$-th element of $\pi$, thus obtaining a permutation
covered by $\pi$ in the $a$-vincular pattern poset. We can repeat
this operation until the first two blocks of $\sigma$ become
adjacent, thus obtaining a chain of coverings in the $a$-vincular
pattern poset from $\pi$ down to a permutation $\overline{\pi}$.
Now observe that, starting from $\overline{\pi}$, we can certainly
remove its leftmost element not belonging to the selected
occurrence of $\sigma$ to obtain a covering, since it is the
$(\alpha_1 +\alpha_2 +1)$-th element of $\pi$ (and $a$ has a 1 in
position $\alpha_1 +\alpha_2$). This argument can be repeated, by
simply observing that, each time we remove the leftmost element of
a permutation of the chain not belonging to the highlighted
occurrence of $\sigma$, we create a new covering in the
$a$-vincular pattern poset. We can thus complete the chain of
coverings from $\tau$ to $\pi$ to a chain of covering from $\tau$
to $\sigma$, which is what we need to conclude that $\sigma \leq_a
\tau$.\cvd

However, also in this special case, the reverse implication does
not hold, as the following example clarifies.

\medskip

\emph{Example.}\quad Let $a=(0,1,0,0,\ldots )$ (that is, the
vector whose unique nonzero component is the second one). If
$\sigma =123$ and $\tau =51423$, then $123\prec_a 4123\prec_a
51423$, and so $\sigma \leq_a \tau$; on the other hand, $\sigma
\notin_a \tau$, since in the unique occurrence of $\sigma$ in
$\tau$ the first two elements are not consecutive.

\medskip

The reason why the above counterexample holds is essentially that
in $a$ there is a 1 preceded by a 0. If this does not happen, we
are able to prove the following proposition.

\begin{prop} Suppose that the vector $a$ has the first $k$ components equal
to 1 and all the remaining components equal to 0. Then $\sigma
\in_a \tau$ if and only if $\sigma \leq_a \tau$.
\end{prop}

\emph{Proof.}\quad We only need to prove that $\sigma \leq_a \tau$
implies that $\sigma \in_a \tau$. Observe that, if the length of
$\sigma$ is $\leq k+1$, then the thesis becomes trivial. So, from
now on, we will tacitly assume that $\sigma$ is sufficiently long.
We proceed by induction on the difference between the length of
$\tau$ and the length of $\sigma$. If $\sigma$ and $\tau$ have the
same length, then clearly $\sigma =\tau$. Moreover, if the length
of $\tau$ is one more than the length of $\sigma$, then the
definitions of $\in_a$ and $\leq_a$ coincide. Now suppose that the
assertion holds for pairs of permutations whose lengths differ by
less than $n$, and let $\sigma$ and $\tau$ be permutations whose
lengths differ by exactly $n$. Since $\sigma \leq_a \tau$, there
is a chain of coverings $\tau =\rho_0 \succ_a \rho_1 \succ_a
\cdots \succ_a \rho_{n-1}\succ_a \rho_n =\sigma$. By inductive
hypothesis, $\rho_{n-1}\in_a \tau$. Due to our assumptions on $a$,
this means that there is an occurrence of $\rho_{n-1}$ in $\tau$
all of whose elements must appear consecutively except at most the
first $k+1$. Since $\sigma \prec_a \rho_{n-1}$, $\sigma$ is
obtained from $\rho_{n-1}$ by removing either one of the first
$k+1$ elements or the last one. In all cases, it is easy to check
that the occurrence of $\sigma$ in $\tau$ thus resulting from the
selected occurrence of $\rho_{n-1} $ in $\tau$ is made of
consecutive elements of $\tau$, except at most the first $k+1$.
This is equivalent to say that $\sigma \in_a \tau$, as
desired.\cvd

\section{The quasi-consecutive pattern poset}\label{quasi}

In the rest of the paper we will deal with a special case, which
turns out to be particularly manageable, due to its closeness to
the consecutive case.

\medskip

Let $A$ be the infinite lower triangular $(0,1)$-matrix for which
an entry is equal to 1 if and only if it belongs to the first
column. So the first lines of $A$ are as follows:

\begin{displaymath}
A=\left(
\begin{array}{ccccc}
  1 & 0 & 0 & 0 & \cdots \\
  1 & 0 & 0 & 0 & \cdots \\
  1 & 0 & 0 & 0 & \cdots \\
  1 & 0 & 0 & 0 & \cdots \\
  \vdots & \vdots & \vdots & \vdots & \ddots \\
\end{array}
\right) .
\end{displaymath}

Since $A$ has constant columns, according to the notations
introduced in the previous section, it corresponds to the infinite
vector $a=(1,0,0,0,\ldots )$ having 1 only in the first position.
Therefore an $a$-vincular pattern $\sigma$ is interpreted as a
dashed permutation having only one dash, which is placed between
the first two elements of $\sigma$. Due to the last result of the
previous section, in this special case, we have that $\sigma
\leq_a \tau$ if and only if $\sigma \in_a \tau$. In other words,
we can assert that $\sigma \leq_a \tau$ when there is an
occurrence of $\sigma$ in $\tau$ whose first two elements are
possibly not consecutive, whereas all the remaining elements has
to appear consecutively in $\tau$. For this reason, we will say
that $\sigma$ is a \emph{quasi-consecutive pattern} of $\tau$, and
the resulting poset will be called the \emph{quasi-consecutive
pattern poset}. As an example, the permutation $432516$ contains
the quasi-consecutive pattern $231$, since the entries $351$ show
an occurrence of the vincular pattern $2-31$. When a permutation
$\tau$ does not contain the quasi-consecutive pattern $\sigma$, we
say that $\tau$ \emph{avoids} the quasi-consecutive pattern
$\sigma$. For instance, the preceding permutation $432516$ avoids
the quasi-consecutive pattern $123$. Observe that other instances
of vincular patterns with ``ground permutation" $123$ appear in
$432516$ (for example, there are occurrences of both $12-3$ and
$1-2-3$). To have an idea of how intervals look like in the
quasi-consecutive pattern poset, we refer the reader to subsequent
figures.

\medskip

The quasi-consecutive pattern poset has many similarities with the
consecutive pattern poset studied in \cite{BFS}. Clearly it has a
slightly more complicated structure, which would be interesting to
investigate in detail. A first structural property of such a poset
is recorded in the next proposition.

\begin{prop} For any $\tau \in S$, $\tau$ covers at most three
permutations, which are obtained by removing either the first, the
second or the last entry of $\tau$ (and suitably rearranging the
remaining ones).
\end{prop}

\emph{Proof.}\quad If $\tau$ has length $n$, a permutation $\rho$
of length $n-1$ covered by $\tau$ can appear as an occurrence
either consecutive or not. In the former case, $\rho$ has to
appear either as a prefix or a suffix of $\tau$; thus $\rho$ is
obtained by removing either the first or the last entry of $\tau$.
In the latter case, the only possibility is that the first two
entries of $\rho$ are not consecutive in $\tau$; thus $\rho$ is
obtained from $\tau$ by removing the second entry.\cvd

The property proved in the above proposition gives an important
feature of the quasi-consecutive poset, which will be very useful
in the next section.

\section{On the M\"obius function of the quasi-consecutive pattern poset}\label{mobius_quasi}

In this section, which contains the main results of our paper, we
address the problem of the computation of the M\"obius function of
the quasi-consecutive pattern poset. Specifically, given $\sigma
,\tau \in S$, with $\sigma\leq \tau$, we want to compute $\mu
(\sigma ,\tau )$. We will completely solve the case of a single
occurrence, which turns out to be not a trivial one, and leave the
general case as an open problem.

Throughout this section, we will frequently write permutations with a
dash between the first entry and the second entry. This is done in
order to emphasize the fact that they have to be interpreted as
specific vincular patterns in the quasi-consecutive pattern poset.
Moreover, as we already started to do in the previous paragraph,
the partial order relation of the quasi-consecutive pattern poset
will be denoted simply by $\leq$ (instead of $\leq_a$), as no
other partial orders will appear from now on.

\medskip

We start by giving some partial results, depending on whether
$\tau$ covers precisely $i$ permutations, for $i=1,2$.
Recall that a \emph{monotone permutation} is a permutation which is
either increasing or decreasing (so that, for any $n$, the monotone
permutations are simply $12\cdots n$ and $n(n-1)\cdots 1$).

\begin{prop} If $\tau$ covers precisely one permutation, then
$[\sigma ,\tau ]$ is a chain, for any $\sigma \leq \tau$.
\end{prop}

\emph{Proof.}\quad If $\tau$ covers precisely one permutation,
then the same permutation is obtained by removing any of the three
allowed elements of $\tau$. In particular, if removing the first
or the last entry of $\tau$ results in the same permutation,
necessarily $\tau$ has to be monotone. So $\sigma$ has to be
monotone too, and in this case $[\sigma ,\tau ]$ is clearly a
chain.\cvd

Thus, if $\tau$ covers precisely one permutation, then $\mst =1$
when $\sigma =\tau$, $\mst =-1$ when $\tau$ covers $\sigma$, and
$\mst =0$ in all the other cases.

%
%

\begin{prop} Let $\tau =a_1 -a_2 \cdots a_n$. If $\tau$ covers precisely two permutations,
then there are two distinct possibilities:
\begin{itemize}

\item $a_1$ and $a_2$ are consecutive integers;

\item either $\tau =1n(n-1)\cdots 32$ or $\tau =n12\cdots (n-1)$.

\end{itemize}
\end{prop}

\emph{Proof.}\quad If $\tau$ covers precisely two permutations,
then there exist two elements in the set $\{ a_1 ,a_2 ,a_n \}$
such that removing either of them from $\tau$ gives the same
permutation. It cannot happen that removing either $a_1$ or $a_n$
yields the same permutation, since otherwise $\tau$ would be
monotone, and so it would cover only one permutation, which is not
the case. If removing either $a_1$ or $a_2$ yields the same
permutation, then each of the remaining entries of $\tau$ has to
be either bigger or smaller than both $a_1$ and $a_2$. Thus $a_1$
and $a_2$ must be consecutive integers. Finally, if removing
either $a_2$ or $a_n$ yields the same permutation, then the
subpermutation determined by the elements $a_2,\ldots ,a_n$ is
monotone and $a_1$ must be either smaller or bigger than such
elements. Since $\tau$ itself cannot be monotone, then $a_1$ is
either 1 or $n$, and so either $\tau =1n(n-1)\cdots 32$ or $\tau
=n12\cdots (n-1)$.\cvd

Observe that, if $\tau$ covers precisely two permutations and
$a_1$ and $a_2$ are not consecutive integers, then we can
determine the structure of the interval $\ist$, and so its
M\"obius function. Indeed, suppose that any occurrence of $\sigma$
in $\tau$ does not involve $a_1$. Then necessarily $\sigma$ is
monotone. Thus a permutation in $\ist$ is uniquely determined by a
(possibly empty) set of entries to remove from the end of $\tau$
and by choosing if $a_1$ has to be removed or not. For these
reasons, it is not difficult to realize that $\ist \simeq
\mathbf{2}\times \mathbf{n}$, where $\mathbf{i}$ denotes the chain
of the first $i$ positive integers. Therefore, in this situation,
the computation of the M\"obius function can be easily carried
out:
\begin{displaymath}
\mst =\left\{ \begin{array}{lll}
1,&\qquad & \textnormal{if $\sigma$ has length $n$;}\\
-1,&\qquad& \textnormal{if $\sigma$ has length $n-1$;}\\
1,&\qquad& \textnormal{if $\sigma$ has length $n-2$;}\\
0,&\qquad& \textnormal{otherwise.}
\end{array}\right.
\end{displaymath}

Notice that $\sigma =1$ does not fit into the previous case,
nevertheless using a completely analogous argument we get to the
same formula for $\mst$.

On the other hand, if there is an occurrence of $\sigma$ which
involves $a_1$, then necessarily every occurrence of $\sigma$
involves $a_1$. In this case, it is easily seen that either
$\sigma$ is not monotone or $\sigma$ has length 2. Thus $\ist$ is
a finite chain, whose M\"obius function is trivial to compute.

\subsection{The case of one occurrence}

Let $\sigma \leq \tau =a_1 -a_2 \cdots a_n$ and suppose that
$\sigma$ occurs precisely once in $\tau$. We will be concerned
with the computation of $\mst$ in this particular case, which we
will be able to solve completely.

We will distinguish several cases, depending on the fact that the
single occurrence of $\sigma$ in $\tau$ does or does not involve
the three elements $a_1 ,a_2$ and $a_n$.

\medskip

A first obvious consideration is that $\sigma$ cannot involve all
of the three above elements, unless $\sigma =\tau$, and in this
last case of course $\mst =1$.

The second possibility is that two of the three elements $a_1 ,a_2
,a_n$ are involved in the unique occurrence of $\sigma$ in $\tau$.
There are of course three distinct cases; however in any of them
the interval $\ist$ is a chain. For instance, if $\sigma$ involves
$a_1$ and $a_n$ (and not $a_2$), then $\sigma$ is (isomorphic to)
$a_1 -a_k a_{k+1}\cdots a_n$, for some $k>2$, and it is immediate
to see that any $\rho \in \ist$ is obtained by starting from
$\tau$ and repeatedly removing the second element of the resulting
permutation until we get to $\rho$. This is of course the only
possible way to remove elements from $\tau$ and remain inside
$\ist$ ($a_1$ and $a_n$ cannot be removed, being part of the
unique occurrence of $\sigma$ in $\tau$). This means that the
interval $\ist$ is a chain. The remaining cases can be dealt with
in a completely analogous way. Thus, we can conclude that, in all
these cases, $\mst =0$, unless $\sigma$ has length $n-1$, in which
case $\mst =-1$.

\medskip

We next examine the case in which only one among $a_1 ,a_2$ and
$a_n$ is involved in the unique occurrence of $\sigma$ in $\tau$.
There are three distinct cases to consider.

\begin{prop}\label{a1} Suppose that $a_1$ occurs in $\sigma$ (whereas $a_2$ and $a_n$
do not). Then $\mst =0$, unless $[\sigma ,\tau ]$ has rank 2 (in
which case $\mst =1$).
\end{prop}

\emph{Proof.}\quad In this case $\sigma$ is (isomorphic to) $a_1
-a_k a_{k+1}\cdots a_{k+h}$, for some $k,h$ with $2<k<n$ and
$0\leq h<n-k$. If $\rho \in \ist$, then $\rho$ is (isomorphic to)
$a_1 -a_i a_{i+1}\cdots a_{k+h+j}$, for suitable $i$ and $j$, and
so it is uniquely determined by two intervals of elements $\{ a_2
,\ldots a_{i-1} \}$ and $\{ a_{k+h+j+1},\ldots a_n \}$ to be
removed from $\tau$ in a well-specified order. Therefore $\ist$ is
a \emph{grid}, i.e. it is isomorphic to a product of two chains,
whose lengths (which of course depend on the values of $i$ and
$j$) are $\geq 1$ (this is due to the fact that $\sigma$ does not
involve $a_2$ and $a_n$). Thus $\mst =0$, unless the two chains
both have length 1. Clearly, in the latter case $\mst =1$.\cvd

\begin{prop} Suppose that $a_2$ occurs in $\sigma$ (whereas $a_1$ and $a_n$
do not). Then $\mst =0$, unless $[\sigma ,\tau ]$ has rank 2 (in
which case $\mst =1$).
\end{prop}

\emph{Proof.}\quad The argument is completely analogous to the one
used in the above proposition. The only difference here is that
the roles of $a_1$ and $a_2$ have to be swapped.\cvd

Notice that the two cases considered so far are essentially
equivalent to the case of a single occurrence in the consecutive
pattern poset. Indeed, using the notation introduced in
\cite{BFS}, it is not too difficult to realize that $\ist$ and
$[\psigma ,\ptau ]$ are order isomorphic.

\medskip

The last case is by far the most challenging one.

\begin{prop}\label{consec_end} Suppose that $a_n$ occurs in $\sigma$ (whereas $a_1$ and $a_2$
do not). Then $\mst =0$, unless $\ist$ has rank 2. In this last
case, if $a_1$ and $a_2$ are consecutive integers, then $\mst =0$,
otherwise $\mst =1$.
\end{prop}

\emph{Proof.}\quad It is convenient to distinguish two cases,
depending on whether $\sigma$ appears as a consecutive pattern in
$\tau$ or not.

If $\sigma$ is not a consecutive pattern of $\tau$, then it is not
difficult to realize that $\ist$ has only one atom, which can be
obtained in the following way: take the unique occurrence of
$\sigma$ in $\tau$ and add to it the element of $\tau$ immediately
to the left of its consecutive part. Observe that, in this case,
if $\sigma$ has length $m$, then $n-m\geq 3$, and $\mst =0$.

\begin{figure}[h!]
\begin{center}
    \includegraphics[scale=.5]{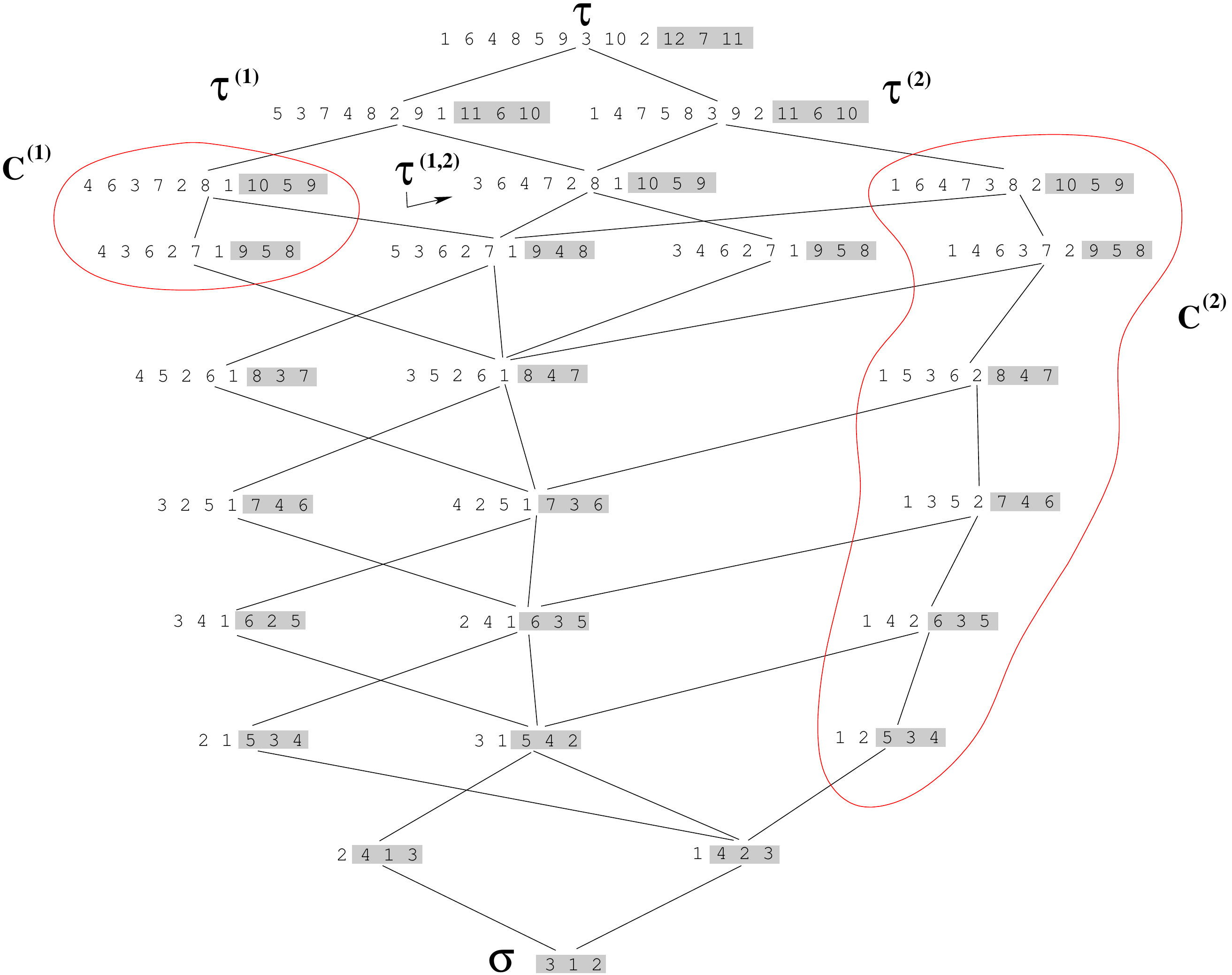}
\end{center}
\caption{A single occurrence of a pattern at the end of a
permutation.}
 \label{prop_4_5}
\end{figure}

If $\sigma$ is a consecutive pattern of $\tau$, then $\sigma$ is
(isomorphic to) $a_k a_{k+1}\cdots a_n$, where $2<k\leq n$. Figure
\ref{prop_4_5} shows an instance of this situation. Denote with
$\tauone$ the permutation obtained from $\tau$ by removing $a_1$
(and of course suitably renaming the remaining elements).
Similarly, denote with $\tautwo$ the permutation obtained from
$\tau$ by removing $a_2$ and with $\tauonetwo$ the permutation
obtained from $\tau$ by removing both $a_1$ and $a_2$. Finally,
let
\begin{eqnarray*}
\Cone &=&\{ \rho \in \ist \; |\, \rho <\tauone ,\rho \nleq
\tauonetwo \} ,\\
\Ctwo &=&\{ \rho \in \ist \; |\, \rho <\tautwo ,\rho \nleq
\tauonetwo \} .
\end{eqnarray*}

We start by observing that, in this case, $\tauone$ and $\tautwo$
are the only coatoms of $\ist$. So, if $\tauone =\tautwo$ (i.e.,
$a_1$ and $a_2$ are consecutive integers), then clearly $\mst =0$.
Thus, in what follows we will suppose that $\tauone \neq \tautwo$.

We next show that $\Cone$ and $\Ctwo$ are chains. Indeed, any
$\rho \in \Ctwo$ can be obtained from $\tautwo$ by repeatedly
removing the second element of the resulting permutation, so
$\rho$ is uniquely determined by the set of consecutive elements
of $\tau$ which have been removed. An analogous argument shows
that also $\Cone$ is a chain.

Moreover, we have that $\Cone \cap \Ctwo =\emptyset$. Indeed, if
we had $\Cone \cap \Ctwo \neq \emptyset$, then any $\rho \in \Cone
\cap \Ctwo$ would be order isomorphic to both $a_1 a_r
a_{r+1}\cdots a_n$ and $a_2 a_r a_{r+1}\cdots a_n$, for a suitable
$r$. We wish to show that, in this situation, $a_1$ and $a_2$ have
to be consecutive integers. Indeed, suppose that $a_1$ and $a_2$
are not consecutive and, w.l.o.g., that $a_1 <a_2$. First of all,
for all $i$ such that $r\leq i\leq n$, it cannot be $a_1 <a_i
<a_2$, otherwise $a_1 a_r a_{r+1}\cdots a_n$ and $a_2 a_r
a_{r+1}\cdots a_n$ would not be order isomorphic. So there should
exist an element $a$ such that $a_1 <a<a_2$, which appears before
$a_r$ in $\tau$. But then $a a_r a_{r+1}\cdots a_n$ would be order
isomorphic to $\rho$, which is not possible, since otherwise $\rho
\leq \tauonetwo$. Therefore $a_1$ and $a_2$ have to be consecutive
integers, which is not true (remember that we are supposing that
$\tauone \neq \tautwo$).
We can thus conclude that $\Cone \cap \Ctwo =\emptyset$, as
desired.

Our next goal is to prove that, for all $\rho \in \Cone \cup
\Ctwo$, $\mu (\rho ,\tau )=0$. In fact, if $\alpha \in [\rho ,\tau
]$, then $\alpha \nleq \tauonetwo$, since otherwise we have $\rho
\leq \tauonetwo$, which is not true. So, if for instance $\rho \in
\Cone$, then we have $\alpha \in \Cone$ as well. Therefore $[\rho
,\tau ]$ is a chain having at least three elements (choose $\alpha =\tauone$),
 whence $\mu (\rho ,\tau )=0$.

Finally, we are now in a position to prove that, for all $\rho
<\tauonetwo$, $\mu (\rho ,\tau )=0$. Indeed, if $\rho$ is covered
by $\tauonetwo$, then necessarily $[\rho ,\tau ]=\{ \rho
,\tauonetwo ,\tauone ,\tautwo ,\tau \}$, and it is immediate to
see that $\mu (\rho ,\tau )=0$. Instead, if $\rho \leq \tauonetwo$
is not covered by $\tauonetwo$, then $[\rho ,\tau ]$ contains the
same five permutations listed above as well as a set $X$ of
permutations less than $\tauonetwo$ (but ``closer" than $\rho$ to
$\tauonetwo$) and (possibly) a set $Y$ of permutations contained
in $\Cone \cup \Ctwo$. Using an inductive argument (on the
distance from $\tauonetwo$), we can show that, for all
permutations $\alpha \in X$, $\mu (\alpha ,\tau )=0$; moreover, we
already know (from the previous paragraph) that, for all $\alpha
\in Y$, $\mu (\alpha ,\tau )=0$. Thus the only $\alpha \in [\rho
,\tau]$ such that $\mu (\alpha ,\tau )\neq 0$ are $\tauonetwo
,\tauone, \tautwo , \tau$, and so we can immediately conclude that
$\mu (\rho ,\tau )=0$.

Since of course $\sigma \leq \tauonetwo$, if the rank of $\ist$ is
greater than 2, then $\sigma < \tauonetwo$ and so $\mst =0$.
Otherwise, if $a_1$ and $a_2$ are consecutive integers, then
$\tauone =\tautwo$, and $\ist$ is a 3-elements chain, so that
$\mst =0$. If instead $a_1$ and $a_2$ are not consecutive, then
$\ist$ is the product of two 2-elements chains, whence $\mst =1$.
\cvd

The last case to be considered is when none of $a_1, a_2$ and
$a_n$ belongs to the unique occurrence of $\sigma$ in $\tau$. It
will be convenient to distinguish two cases, depending on whether
the occurrence of $\sigma$ is consecutive or not.

\begin{prop}\label{none_nonconsecutive} If $\sigma$ does not occur consecutively in $\tau$,
and does not involve $a_1 ,a_2$ and $a_n$, then $\mst =0$.
\end{prop}

\emph{Proof.}\quad Denote with $\pi _1$ the permutation isomorphic
to the smallest consecutive pattern in $\tau$ containing $\sigma$
and with $\pi _2$ the permutation isomorphic to the smallest
suffix of $\tau$ containing $\sigma$.

We start by observing that, for all $\rho \in \ist$, we have $\rho
\leq \pi_2$ or $\rho \geq \pi_1$. Indeed, if $\rho \nleq \pi_2$,
then the leftmost element of an occurrence of $\rho$ in $\tau$
must correspond to an entry of $\tau$ which appears on the left of
$\sigma$. Thus necessarily all the other elements of the
occurrence of $\rho$ have to appear consecutively in $\tau$ and to
contain $\sigma$, and this implies that $\rho \geq \pi_1$.

\begin{figure}[h!]
\begin{center}
    \includegraphics[scale=.4]{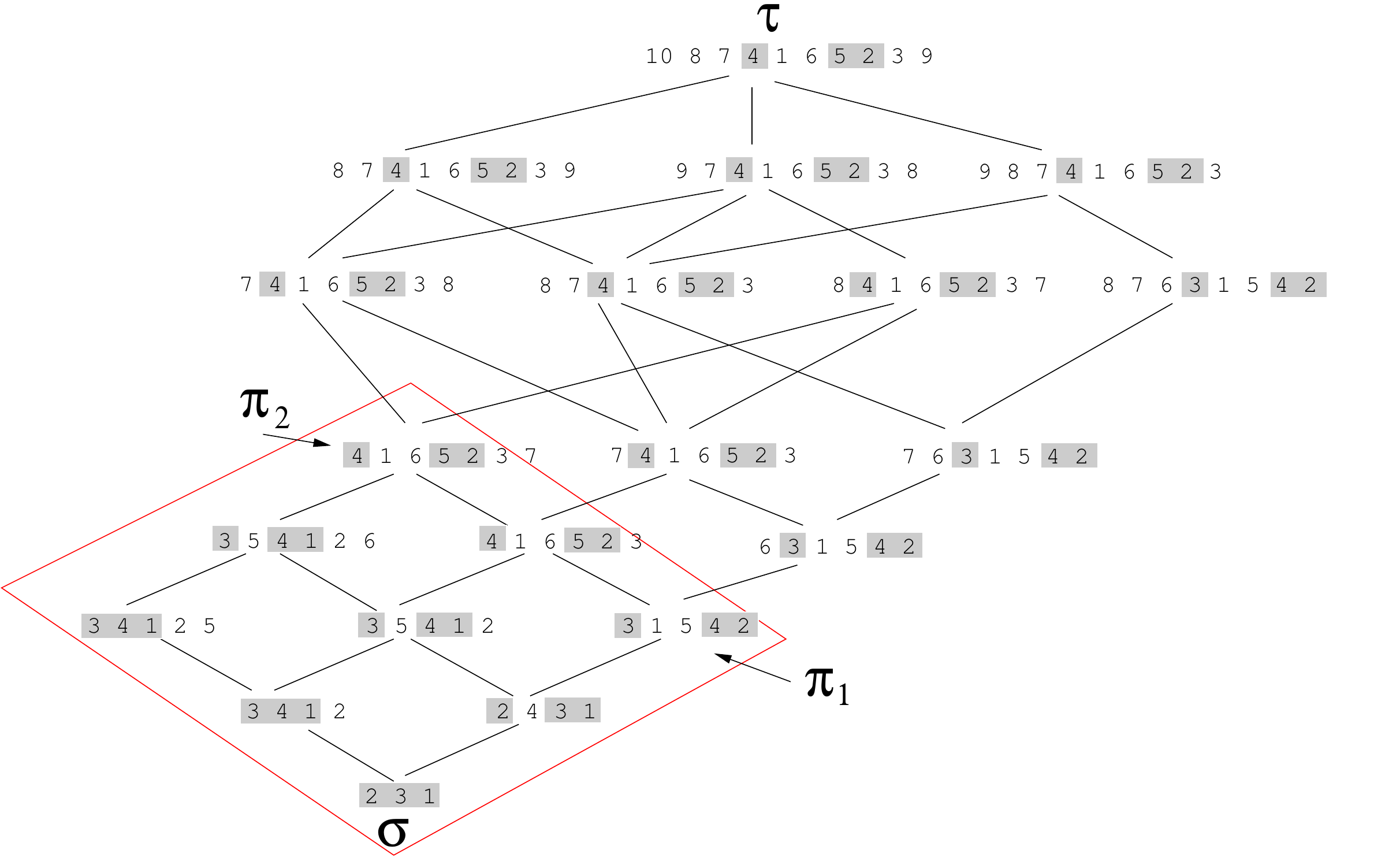}
\end{center}
\caption{
The unique occurrence of $\sigma$ in each permutation of $\ist$ is
highlighted.}
 \label{particular}
\end{figure}

We next give more precise information on the structure of the
interval $\ist$ (see also Figure \ref{particular}).

\begin{itemize}

\item $[\pi_1 ,\pi_2 ]$ is a chain: this follows from the fact
that $\pi_1$ occurs only once in $\pi_2$, and the first two
elements of $\pi_2$ belong to such an occurrence.

\item $[\sigma ,\pi_1 ]$ is a chain: this is analogous to the
above one, since $\sigma$ occurs only once in $\pi_1$ and the
first and last elements belong to such an occurrence.

\item $[\sigma ,\pi_2 ]$ is a product of two nonempty chains: this
follows from the proof of Proposition \ref{a1}.

\end{itemize}

%

Therefore we have
$$
\mst =-\sum_{\sigma \leq \rho <\tau}\mu (\sigma ,\rho
)=-\sum_{\rho \in [\sigma ,\pi_2 ]}\mu (\sigma ,\rho )-\sum_{\rho
\notin [\sigma ,\pi_2 ]}\mu (\sigma ,\rho ).
$$

It is immediate to see that, in the last expression, the first sum
is 0. As for the second sum, we observe that each $\rho \geq
\pi_1$ (and so in particular each $\rho \notin [\sigma ,\pi_2 ]$)
lies either above all the four elements of $[\sigma ,\pi_2 ]$
having nonzero M\"obius function, or above two of them, namely
$\sigma$ and one of the two atoms of $\ist$. In both cases, an
inductive argument (on the difference between the length of $\rho$
and the length of $\pi_1$) shows that, for all $\rho \notin
[\sigma ,\pi_2 ]$, $\mu (\sigma ,\rho )=0$. From all the above
considerations it then follows that $\mst =0$.\cvd

\begin{prop}\label{ultima} If $\sigma$ occurs consecutively in $\tau$,
and does not involve $a_1 ,a_2$ and $a_n$, then $\mst =0$, unless
$\ist$ has rank 3, the unique occurrence of $\sigma$ in $\tau$
contains $a_{n-1}$ and $\tau$ covers three elements in $\ist$. In
such a case, we have $\mst =-1$.
\end{prop}

\emph{Proof.}\quad Denote with $\pi$ the permutation isomorphic to
the smallest suffix of $\tau$ containing the unique occurrence of
$\sigma$ and with $\eta$ the permutation isomorphic to the prefix
of $\tau$ of length $n-1$.

Using an argument analogous to that of Proposition
\ref{none_nonconsecutive}, we observe that, for every $\rho \in
\ist$, we have $\rho \leq \eta$ or $\rho \geq \pi$. Indeed, if
$\rho \nleq \eta$, then an occurrence of $\rho$ necessarily
contains both the last element of $\tau$ and the unique occurrence
of $\sigma$, so that $\rho$ contains $\pi$ (i.e. $\rho \geq \pi$).

Since $\pi$ and $\eta$ are clearly incomparable, the two intervals
$[\sigma ,\eta ]$ and $[\pi ,\tau ]$ constitute a partition of
$\ist$. Therefore concerning the M\"obius function of $\ist$ we
have:
$$
\mst =-\sum_{\pi \leq \rho <\tau}\mu (\sigma ,\rho ).
$$

\begin{figure}[h!]
\begin{center}
    \includegraphics[scale=.5]{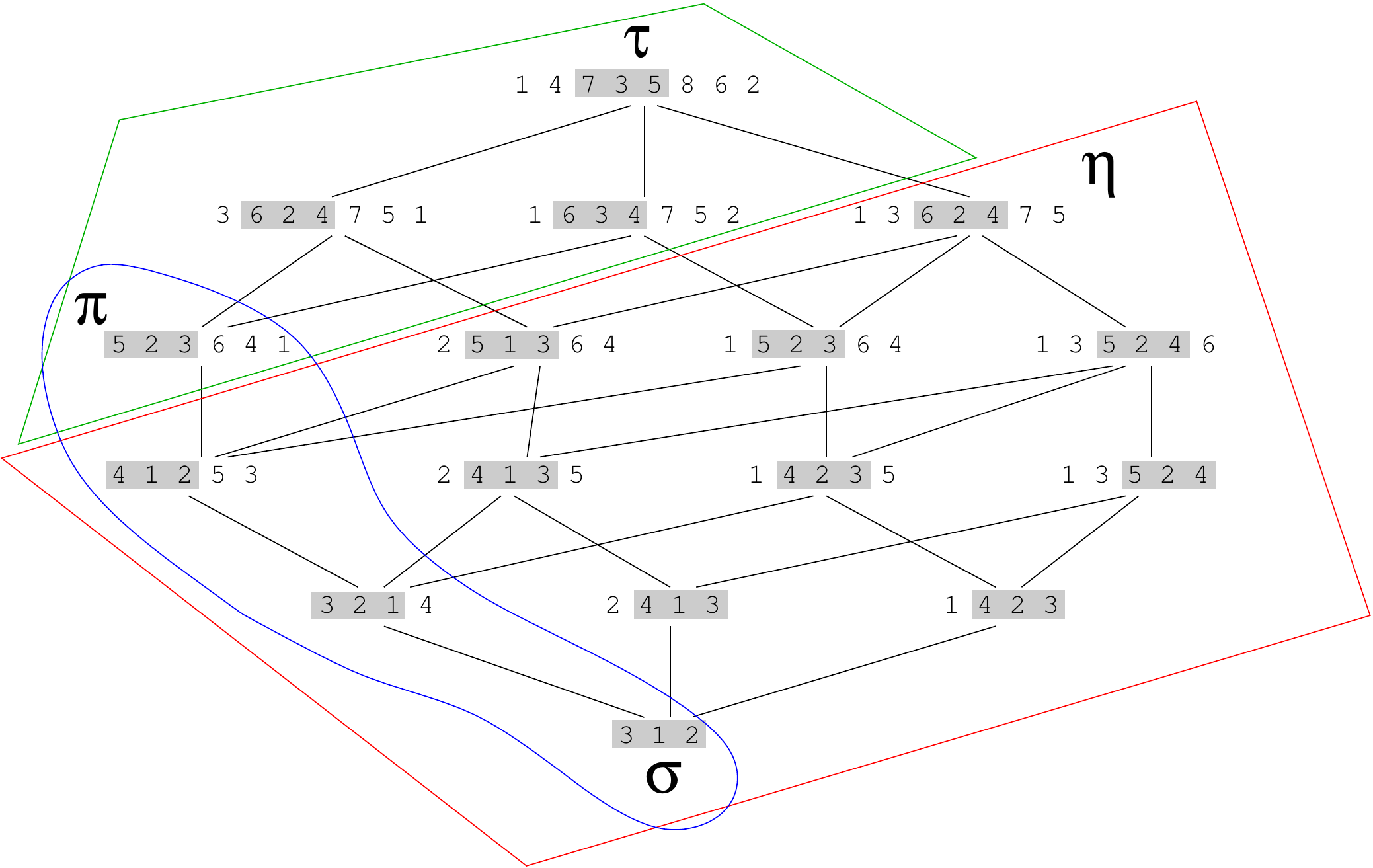}
\end{center}
\caption{
A single occurrence of a consecutive pattern not involving the
three critical entries of a permutation (that is, the first two
and the last one).}
 \label{prop_4_7}
\end{figure}

Suppose first that $[\sigma ,\pi ]$ has rank at least 2. This
situation is illustrated in Figure \ref{prop_4_7}. Observe that,
in this case, $\mu (\sigma ,\pi )=0$, since $\sigma$ is a prefix
of $\pi$, and so $[\sigma ,\pi ]$ is a chain of length at least 2.
Given $\rho \in [\pi ,\tau ]$, suppose now that $\mu (\sigma
,\alpha )=0$, for all $\pi \leq \alpha <\rho$. Notice that $\rho$
is made by an element of $\tau$ followed by a suffix of $\tau$
containing $\sigma$. Thus, by removing the rightmost element of
$\rho$, we obtain a permutation $\beta$ contained in $\eta$
(observe that $\beta$ is a coatom of $[\sigma ,\rho ]$). Moreover,
the reader can immediately see that $\beta$ and $\pi$ are
incomparable and that, for all $\alpha \in [\sigma ,\rho ]$,
either $\alpha \leq \beta$ or $\alpha \geq \pi$. So we have a
partition of $[\sigma ,\rho ]$ and we can conclude that, for all
$\rho \in [\pi ,\tau ]$,
$$
\mu (\sigma ,\rho )=-\sum_{\sigma \leq \alpha <\rho}\mu (\sigma
,\alpha )=-\sum_{\pi \leq \alpha <\rho}\mu (\sigma ,\alpha )=0,
$$
which is enough to assert that $\mst =0$.

Otherwise, i.e. when $[\sigma ,\pi ]$ has rank 1, we observe that
the two intervals $[\sigma ,\eta ]$ and $[\pi ,\tau ]$ both fall
into the scopes of Proposition \ref{consec_end}. Indeed, there is
a unique occurrence of $\sigma$ and $\pi$ in $\eta$ and $\tau$
respectively, in both cases at the end of the permutation. Each
permutation $\gamma \in [\pi ,\tau ]$ covers exactly one
permutation $\gamma'\in[\sigma ,\eta ]$ (just remove the last
entry of $\gamma$, which is its only entry which can be removed in
order to get a permutation lying below $\eta$). This fact will now
be used to describe the structure of $\ist$.

\begin{figure}[h!]
\begin{center}
    \includegraphics[scale=.5]{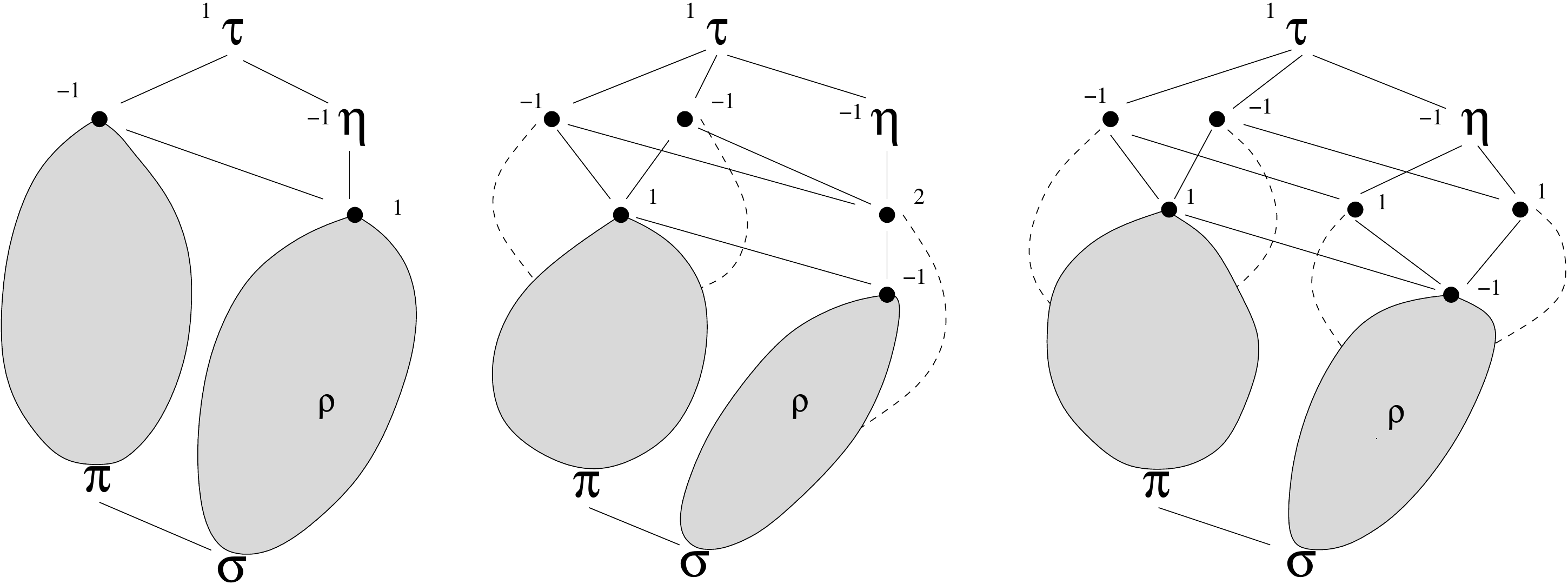}
\end{center}
\caption{The M\"obius function is computed from top to bottom.
Dashed lines denote chains. In each of the three cases, there are
covering relations between elements of the two shaded regions, as
well as between elements of two corresponding chains, under the
map denoted with $'$ in the proof of Proposition \ref{ultima}.}
 \label{prop_4_7_2}
\end{figure}

Concerning the interval $[\pi ,\tau ]$, we have two possibilities:
either $\tau$ covers one element or two elements. In the former
case, necessarily $\eta$ covers only one element in $[\sigma ,\eta
]$; in the latter case, $\eta$ may cover either one or two
elements. We thus have a total of three cases, which are
illustrated in Figure \ref{prop_4_7_2}.

We first observe that, when $\ist$ has rank 3, the second and the
third cases depicted in Figure \ref{prop_4_7_2} correspond exactly
to the exceptional case of the statement of the proposition, and
it is immediate to observe that $\mst =-1$. Instead, in the first
case we have $\mst =0$. So from now on we suppose that the rank of
$\ist$ is at least 4.

If we compute the M\"obius function of $\ist$ from top to bottom,
we can restrict ourselves to compute it in each element of
$[\sigma ,\eta ]$, that is:
$$
\mst =-\sum_{\sigma <\rho \leq \eta}\mu (\rho ,\tau ).
$$

It is now not too difficult to realize that, in each of the three
possible cases, for any $\rho \in [\sigma ,\eta ]$ except that for
the top three levels of the interval (i.e. $\eta$, its coatoms and
the elements covered by such coatoms), the value of $\mu (\rho
,\tau )$ is 0. Indeed, denote with $\Delta$ the set of all
permutations in $[\rho ,\tau ]\setminus \{ \rho \}$ other than
those whose M\"obius function is explicitly indicated in Figure
\ref{prop_4_7_2}. Then, $$\mu (\rho ,\tau )=-\sum_{\delta \in
\Delta}\mu (\delta ,\tau )-\sum_{\xi >\rho \atop \xi\notin
\Delta}\mu (\xi ,\tau ).$$ We can assume inductively that each
summand of the first sum is equal to 0; as far as the second sum
is concerned, the reader is referred again to Figure
\ref{prop_4_7_2} to get convinced that it is equal to 0 in all
cases as well. Thus, in particular, $\mst =0$.\cvd

Summing up all the results obtained so far in this section, we
then have the following theorem, which completely solves the
problem of the computation of the M\"obius function of the
quasi-consecutive pattern poset in the case of one occurrence.

\begin{teor} Suppose that $\sigma$ occurs exactly once in $\tau =a_1 -a_2 \cdots a_n$
as a quasi-consecutive pattern. Then $\mst =0$, unless one of the
following cases hold:
\begin{itemize}

\item $\sigma =\tau$, in which case $\mst =1$;

\item $\sigma$ is covered by $\tau$, in which case $\mst =-1$.

\item $\ist$ has rank 2 and $\sigma$ involves $a_1$ but not $a_2$
and $a_n$, in which case $\mst =1$ (the same holds when $a_1$ and
$a_2$ are swapped).

\item $\ist$ has rank 2, $\sigma$ involves $a_n$ but not $a_1$ and
$a_2$, and $a_1$ and $a_2$ are not consecutive integers, in which
case $\mst =1$.

\item $\ist$ has rank 3, $\sigma$ involves $a_{n-1}$ but not $a_1
,a_2$ and $a_n$, $\sigma$ occurs consecutively in $\tau$ and
$\tau$ covers three elements in $\ist$, in which case $\mst =-1$.

\end{itemize}

\end{teor}




\section{Further work}

The study of vincular pattern posets, which has been initiated in
the present paper, is of course very far from being completed.

From the point of view of vincular pattern posets in general, the
main open problem, already stated in Section \ref{posets}, is
perhaps that of characterizing those matrices $A$ for which
$\sigma \leq_A \tau$ if and only if $\sigma \in_A \tau$.

Concerning the main topic investigated here, namely the M\"obius
function of the quasi-consecutive pattern poset, a lot of work has
still to be done. The case of one occurrence which we have
completely solved in the previous section seems not to be really
representative of the general case. For instance, the absolute
value of $\mu$ can be different from 0 and 1 (a simple example is
given by the interval $[12,2413]$, whose M\"obius function equals
2), and there is computational evidence \cite{St} that $|\mu
(\sigma ,\tau )|$ is actually unbounded. Despite its closeness
with the consecutive case, this fact shows that the
quasi-consecutive case can sometimes be very similar to the
classical (unrestricted) case. Another conjecture suggested by
\cite{St} is the following: if $\tau$ is the direct sum of some
copies of $\sigma$, then $\mst =1$.

\medskip

{\footnotesize {\bf Acknowledgment.}\quad L.F. wishes to thank
Einar Steingr\'imsson for his kind hospitality at the University
of Strathclyde in November 2012. Some of the results and
conjectures of the present work have their roots in several
stimulating discussions during that visit.}

\end{document}